\documentclass[12pt]{article}

\usepackage{amssymb}

\usepackage[english,francais]{babel}

\newtheorem{e-proposition}[theorem]{Proposition}

\newtheorem{e-definition}[theorem]{Definition\rm}

\title{An extension to the Wiener space of the arbitrary functions principle}

\author{Nicolas Bouleau\\ ENPC, ParisTech}
\date{---}
\begin{document}
\maketitle

\selectlanguage{francais}



\selectlanguage{english}

\begin{abstract}


The arbitrary functions principle says that the fractional part of $nX$ converges stably to an independent  random variable uniformly distributed on the unit interval, as soon as the random variable $X$ possesses a density or a characteristic function vanishing at infinity. We prove a similar property for random variables defined on the Wiener space when the stochastic measure $dB_s$ is crumpled on itself.
 {\it To cite this article: N. Bouleau, C. R.
Acad. Sci. Paris, Ser. I 343, (2006), 329-332.}

\vskip 0.5\baselineskip

\selectlanguage{francais}


\noindent{\bf R\'esum\'e}

\vskip 0.5\baselineskip

\noindent
Le principe des fonctions arbitraires dit que la partie fractionnaire de $nX$ converge stablement vers une variable al\'eatoire ind\'ependante uniform\'ement r\'epartie sur $[0,1]$ d\`es que $X$ a une densit\'e ou seulement une fonction caract\'eristique tendant vers z\'ero \`a l'infini. Nous \'etablissons une propri\'et\'e analogue pour des variables al\'eatoires d\'efinies sur l'espace du mouvement brownien par repliement de la mesure stochastique $dB_s$ sur elle-m\^eme.
 {\it Pour citer cet article~: N. Bouleau, C. R.
Acad. Sci. Paris, Ser. I 343, (2006), 329-332.}

\end{abstract}

\selectlanguage{english}


\section{Introduction}
Let us denote $\{x\}$ the fractional part of the real number $x$ and $\stackrel{d}{\Longrightarrow}$  the weak convergence of random variables.  Let $(X,Y)$ be a pair of random variables with values in $\mathbb{R}\times\mathbb{R}^r$, we refer to the following property or its extensions as the arbitrary functions principle:
\begin{equation}(\{nX\},Y)\quad\stackrel{d}{\Longrightarrow}\quad (U,Y)\end{equation}
where $U$ is uniformly distributed on $[0,1]$ independent of $Y$.

This property is satisfied when $X$ has a density or more generally a characteristic function vanishing at infinity. (cf [5] Chap. VIII \S92  and \S93, [2], [4]). It yields an approximation property of $X$ by the random variable $X_n= X-\frac{1}{n}\{nX\}=\frac{[nX]}{n}$ where $[x]$ denotes the entire part of $x$:\\

\noindent{\bf Proposition 1.} {\it Let $X$ be a real random variable with density and  $Y$ a random variable with values in  $\mathbb{R}^r$. Let $X_n=\frac{[nX]}{n}$

a) For all $\varphi\in\mathcal{C}^1\bigcap{\mbox{\rm Lip}}(\mathbb{R})$ and   for all integrable random variable $Z$, 
$$(n(\varphi(X_n)-\varphi(X)), Y)\quad\stackrel{d}{\Longrightarrow}\quad(-U\varphi^\prime(X),Y)$$
$$n^2\mathbb{E}[(\varphi(X_n)-\varphi(X))^2Z]\quad\rightarrow\quad\frac{1}{3}\mathbb{E}[\varphi^{\prime 2}(X)Z]$$
where $U$ is uniformly distributed on $[0,1]$ independent of $(X,Y)$.

\indent b) $\forall\psi\in L^1([0,1])$
$$(\psi(n(X_n-X)),Y)\quad\stackrel{d}{\Longrightarrow}\quad(\psi(-U),Y)$$
under any probability measure $\tilde{\mathbb{P}}\ll\mathbb{P}$.}

We extend such results to random variables defined on the Wiener space. 
\section{Periodic isometries.}

Let $(B_t)$ be a standard $d$-dimensional Brownian motion and let $m$ be the Wiener measure, law of $B$. Let $t\mapsto M_t$ be a bounded deterministic measurable map, periodic with unit period, into the space of orthogonal $d\times d$-matrices such that $\int_0^1M_s ds=0$ (e.g. a rotation in $\mathbb{R}^d$ of angle $2\pi t$).  The transform $B_t\mapsto\int_0^tM_sdB_s$ defines an isometric endomorphism in $L^p(m), 1\leq p\leq \infty$. Let be $M_n(s)=M(ns)$ and $T_n=T_{M_n}$. The transposed of the matrix $N$ is denoted $N^\ast$.\\

\noindent{\bf Proposition 2.} {\it Let be $X\in L^1(m)$. Let $\tilde{m}$ be a probability measure  absolutely continuous w.r. to $m$. Under $\tilde{m}$ we have
$$(T_n(X),B)\quad\stackrel{d}{\Longrightarrow}\quad(X(w),B).$$
The weak convergence acts on $\mathbb{R}\times\mathcal{C}([0,1])$ and $X(w)$  denotes a random variable with the same law as $X$ had under $m$ function of a Brownian motion $W$ independent of $B$.}\\

\noindent{\bf Proof.} a) If $X=\exp\{i\int_0^1\xi.dB+\frac{1}{2}\int_0^1|\xi|^2ds\}$ for some element $\xi\in L^2([0,1],\mathbb{R}^d)$, we have
$T_n(X)=\exp\{i\int_0^1\xi^\ast_sM_n(s)dB_s+\frac{1}{2}\int_0^1|\xi|^2ds\}.$

Putting $Z^n_t=\int_0^t\xi^\ast_sM_n(s)dB_s$  gives
$\langle Z^n,Z^n\rangle_t=\int_0^t\xi^\ast_sM_n(s)M_n^\ast(s)\xi_sds=\int_0^t|\xi|^2(s)ds$ which is a continuous function. Now by proposition 1, 
$\int_0^t\xi^\ast_sM_n(s)ds\rightarrow\int_0^t\xi^\ast_sds\int_0^1M_n(s)ds=0.$
which implies by Ascoli theorem $sup_t|\int_0^t\xi^\ast_sM_n(s)ds|\rightarrow0.$ The argument of H. Rootz\'en  [6] applies and yields
$(\int_0^.\xi^\ast M_ndB,B)\stackrel{d}{\Longrightarrow}(\int_0^.\xi.dW,B)$ giving the result in this case by continuity of the exponential function.

b) When $X\in L^1(m)$, we approximate $X$ by $X_k$ linear combination of exponentials of the preceding type and consider the caracteristic functions. The inequality $$|\mathbb{E}[e^{iuT_n(X)}e^{i\int h.dB}-\mathbb{E}[e^{iuT_n(X_k)}e^{i\int h.dB}]|
\leq |u|\mathbb{E}|T_n(X)-T_n(X_k)|=|u|\;\|X-X_k\|_{L^1}$$
 gives the result.
 
 c) This extends to the case $\tilde{m}\ll m$ by the properties of stable convergence.\hfill$\diamond$

\section{Approximation of the Ornstein-Uhlenbeck structure.}
From now on, we assume for simplicity that $(B)$ is one-dimensional. Let $\theta$ be a periodic real function with unit period such that$\int_0^1\theta(s)ds=0$ and $\int_0^1\theta^2(s)ds=1$. We consider the transform $R_n$ of the space $L^2_\mathbb{C}(m)$ defined by its action on the Wiener chaos: 

If $X=\int_{s_1<\cdots<s_k}\hat{f}(s_1,\ldots,s_k)dB_{s_1}\ldots dB_{s_k}$ for $\hat{f}\in L^2_{sym}([0,1]^k,\mathbb{C})$,
$$R_n(X)=\int_{s_1<\cdots<s_k}\hat{f}(s_1,\ldots,s_k)e^{i\frac{1}{n}\theta(ns_1)}dB_{s_1}\ldots e^{i\frac{1}{n}\theta(ns_k)}dB_{s_k}.$$
$R_n$ is an isometry from $L^2_\mathbb{C}(m)$ into itself. From
$n(e^{\frac{i}{n}\sum_{p=1}^k\theta(ns_p)}-1)=i\sum_{p=1}^k\theta(ns_p)\int_0^1 e^{\alpha \frac{i}{n}\sum_p\theta(ns_p)}d\alpha$ it follows that if $X$ belongs to the $k$-th chaos
$$\|n(R_n(X)-X)\|^2_{L^2}\leq k^2\|X\|^2_{L^2}\|\theta\|^2_\infty.$$ In other words, denoting $A$ the Ornstein-Uhlenbeck operator, $X\in \mathcal{D}(A)$ implies
$$\|n(R_n(X)-X)\|_{L^2}\leq 2\|AX\|_{L^2}\|\theta\|_\infty$$ and this leads to\\ 

\noindent{\bf Proposition 3.} {\it If $X\in\mathcal{D}(A)$
$$(-in(R_n(X)-X),B)\quad\stackrel{d}{\Longrightarrow}\quad(X^\#(\omega,w),B)$$
where $W$ is an Brownian motion independent of $B$ and $X^\#=\int_0^1D_sX\,dW_s$.}\\

\noindent{\bf Proof.} {If $X$ belongs to the $k$-th chaos, expanding the exponential by its Taylor series gives 
$$n(R_n(X)-X) =i\int_{s_1<\cdots<s_k}\hat{f}(s_1,\ldots,s_k)\sum_{p=1}^k\theta(ns_p)dB_{s_1}\ldots dB_{s_k}+Q_n$$
with $\|Q_n\|^2\leq \frac{1}{4n}k^2\|\theta\|^2_\infty\|X\|^2$.

Then using that $\int_{s_1<\cdots<s_p<\cdots<s_k}h(s_1,\ldots,s_k)\theta(ns_p)dB_{s_1}\ldots dB_{s_p}\ldots dB_{s_k}$ 

\noindent converges stably to $\int_{s_1<\cdots<s_p<\cdots<s_k}h(s_1,\ldots,s_k)dB_{s_1}\ldots dW_{s_p}\ldots dB_{s_k}$ one gets 
$$-in(R_n(X)-X)\quad\stackrel{s}{\Longrightarrow}\quad
\begin{array}{l}
\int_{t<s_2<\cdots<s_k}\hat{f}(t,s_2,\ldots,s_k)dW_tdB_{s_2}\ldots dB_{s_k}\\
+\int_{s_1<t<\cdots<s_k}\hat{f}(s_1,t,\ldots,s_k)dB_{s_1}dW_t\ldots dB_{s_k}\\
+\cdots\\
+\int_{s_1<\cdots<s_{k-1}<t}\hat{f}(s_1,\ldots,s_{k-1},t)dB_{s_1}\ldots dB_{s_{k-1}}dW_t
\end{array}
$$
which equals $\int D_s(X)dW_s=X^\#$.

The general case in obtained by approximation of $X$ by $X_k$ for the $\mathbb{D}^{2,2}$ norm and the same argument as in the proof of proposition 2 by the caracteristic functions gives the result.\hfill$\diamond$

By the properties of stable convergence, the weak convergence of prop. 3 also holds under $\tilde{m}\ll m$. By similar computations we obtain\\

\noindent{\bf Proposition 4.} {\it  $\forall X\in\mathcal{D}(A)$
$$n^2\mathbb{E}[|R_n(X)-X|^2]\rightarrow2\mathcal{E}[X]$$
where $\mathcal{E}$ is the Dirichlet form associated with the Ornstein-Uhlenbeck  operator.}\\

Following the same lines, it is possible to show that the theoretical  $\overline{A}$ and practical  $\underline{A}$ bias operators (cf. [1])  defined on the algebra $\mathcal{L}\{e^{\int\xi dB}\;;\;\xi\in\mathcal{C}^1\}$ by 
$$\begin{array}{c}
n^2\mathbb{E}[(R_n(X)-X)Y]=<\overline{A}X,Y>_{L^2(m)}\\
n^2\mathbb{E}[(X-R_n(X))R_n(Y)]=<\underline{A}X,Y>_{L^2(m)}
\end{array}
$$ are defined and equal to $A$.\\

\noindent{\it Comment.} The preceding properties are very similar to the results concerning the weak asymptotic error for the resolution of SDEs  by the Euler scheme, involving also an ``extra"-Brownian motion (cf. [3]).

Nevertheless these results do not use the arbitrary functions principle  because a convergence like 
$(n\int_0^.(s-\frac{[ns]}{n})dB_s,B)\stackrel{d}{\Longrightarrow}(\frac{1}{\sqrt{12}}W+\frac{1}{2}B,B)$
  is hidden by a dominating phenomenon
$(\sqrt{n}\int_0^.(B_s-B_{\frac{[ns]}{n}}dB_s,B)$ $\stackrel{d}{\Longrightarrow}(\frac{1}{\sqrt{2}}
\tilde{W},B)$
due to the fact that when a sequence of variables in the second (or higher order) chaos converges stably to a Gaussian variable, this one appears to be independent of the fisrt chaos and therefore of $B$.

The arbitrary functions principle is  slightly different, it is a crumpling of the random orthogonal measure $dB_s$ on itself. This operates even on the first chaos.  Concerning the solution of SDEs by the Euler scheme, it is in force   for SDEs of the form 
$$\left\{
\begin{array}{l}
X^1_t=x^1_0+\int_0^tf^{11}(X^2_s)dB_s+\int_0^tf^{12}(X^1_s,X^2_s)ds\\
X^2_t=x^2_0+\int_0^tf^{22}(X^1_s,X^2_s)ds
\end{array}\right.
$$
where $X^1$ is with values in $\mathbb{R}^{k_1}$,  $X^2$ in $\mathbb{R}^{k_2}$,  $B$ in $\mathbb{R}^d$ and $f^{ij}$ are matrices with suitable dimensions which are encountered for the description of mechanical systems under noisy sollicitations when the noise depend only on the position of the system and the time. In such equations,  integration by parts reduces the stochastic integrals to ordinary integrals and it may be shown  that solved by the Euler scheme they present a weak asymptotic error in  $\frac{1}{n}$ instead of $\frac{1}{\sqrt{n}}$ as usual.

\end{document}